\documentclass[12pt]{amsart}
\usepackage{latexsym,amssymb,amscd}

\def\NZQ{\Bbb}               
\def\NN{{\NZQ N}}

\def\ZZ{{\NZQ Z}}

\def\B'c{{\mathcal{B'}}}
\def\U'c{{\mathcal{U'}}}

\def\opn#1#2{\def#1{\operatorname{#2}}} 
\opn\chara{char}
\opn\length{\ell}
\opn\projdim{proj\,dim}
\opn\injdim{inj\,dim}
\opn\ini{in}
\opn\rank{rank}
\opn\depth{depth}
\opn\sdepth{sdepth}
\opn\height{ht}
\opn\embdim{emb\,dim}
\opn\codim{codim}

\opn\Tr{Tr}
\opn\bigrank{big\,rank}
\opn\superheight{superheight}\opn\lcm{lcm}
\opn\trdeg{tr\,deg}%
\opn\reg{reg}
\opn\lreg{lreg}
\opn\set{set}
\opn\supp{Supp}
\opn\shad{Shad}
\opn\div{div}
\opn\Div{Div}
\opn\cl{cl}
\opn\Cl{Cl}
\opn\Spec{Spec}
\opn\Supp{Supp}
\opn\supp{supp}
\opn\Sing{Sing}
\opn\Ass{Ass}
\opn\Ann{Ann}
\opn\Rad{Rad}
\opn\Soc{Soc}
\opn\Ker{Ker}
\opn\Coker{Coker}
\opn\Im{Im}
\opn\Hom{Hom}
\opn\Tor{Tor}
\opn\Ext{Ext}
\opn\End{End}
\opn\Aut{Aut}
\opn\id{id}

\opn\nat{nat}
\opn\GL{GL}
\opn\SL{SL}
\opn\mod{mod}
\opn\ord{ord}
\opn\aff{aff}
\opn\con{conv}
\opn\relint{relint}
\opn\st{st}
\opn\lk{lk}
\opn\cn{cn}
\opn\core{core}
\opn\vol{vol}
\opn\gr{gr}

\def\pot#1#2{#1[\kern-0.28ex[#2]\kern-0.28ex]}

\opn\dirlim{\underrightarrow{\lim}}
\opn\invlim{\underleftarrow{\lim}}
\def\pnt{{\raise0.5mm\hbox{\large\bf.}}}

\def\Implies{\ifmmode\Longrightarrow \else
     \unskip${}\Longrightarrow{}$\ignorespaces\fi}
\def\implies{\ifmmode\Rightarrow \else
     \unskip${}\Rightarrow{}$\ignorespaces\fi}
\def\iff{\ifmmode\Longleftrightarrow \else
     \unskip${}\Longleftrightarrow{}$\ignorespaces\fi}

\let\:=\colon
\newtheorem{Theorem}{Theorem}[section]
\newtheorem{Lemma}[Theorem]{Lemma}
\newtheorem{Corollary}[Theorem]{Corollary}
\newtheorem{Proposition}[Theorem]{Proposition}
\newtheorem{Remark}[Theorem]{Remark}

\newtheorem{Example}[Theorem]{Example}

\let\epsilon=\varepsilon
\let\phi=\varphi
\let\kappa=\varkappa
\numberwithin{equation}{section}

\title{upper bounds for the Stanley depth }
\author{muhammad Ishaq}
\thanks{The author would like to express his gratitude to ASSMS of GC University Lahore for creating a very appropriate
atmosphere for research work. This research is partially supported by HEC Pakistan}
\address{Muhammad Ishaq, Abdus Salam School of Mathematical Sciences, GC University, Lahore, 68-B New Muslim Town Lahore, Pakistan.}
\email{ishaq\_maths@yahoo.com}

\begin{document}

\maketitle

\begin{abstract}
Let $I\subset J$ be monomial ideals of a polynomial algebra $S$ over a field. Then the Stanley depth of $J/I$ is smaller or equal with the Stanley depth of $\sqrt{J}/\sqrt{I}$. We give also an upper bound for the Stanley depth of the intersection of two primary monomial ideals $Q$, $Q'$, which is reached  if $Q$, $Q'$ are irreducible, ht$(Q+Q')$ is odd and $\sqrt{Q}$, $\sqrt{Q'}$ have no common variable.2000 Mathematics Classification: Primary 13H10, Secondary 13P10, 13C14, 13F20.

\end{abstract}
\section{Introduction}
Let $S=K[x_1,\dots,x_n]$ be a polynomial ring over a field $K$
and $M$ be a finitely generated $\ZZ^n$-graded $S$-module. For a homogenous element $w\in M$ and a subset $Z\subset \{x_1,\dots,x_n\}$, $wK[Z]$ denotes the $K$-subspace of $M$ generated by all homogeneous elements of the form $wv$, where $v$ is a monomial in $K[Z]$. The linear $K$-subspace $wK[Z]$ is called a Stanley space of dimension $|Z|$ if it is a free $K[Z]$-module, where $|Z|$ denotes the number of indeterminates in $Z$. A Stanley decomposition of $M$ is a presentation of the $K$-vector space $M$ as a finite direct sum of Stanley spaces$$\mathcal{D}:M=\bigoplus_{i=1}^su_iK[Z_i]$$ and the Stanley depth of a decomposition $\mathcal{D}$ is  $\sdepth\mathcal{D}=\min\{|Z_i|,i=1,\dots,s\}$. The Stanley depth of $M$ is  $$\sdepth_S(M)=\max\{\sdepth\mathcal{D}:\mathcal{D} \ \text{is a Stanley decomposition of} \ M\},$$ (sometimes we write $\sdepth(M)$ if no confusion is possible). Several properties of the Stanley depth are given in \cite{BK},\cite{HVZ},\cite{R1},\cite{R2} and \cite{SJ}.  Stanley conjectured in \cite{RP} that $\sdepth(M)\geq \hbox{depth}(M)$ (for terminology of commutative algebra we refer to \cite{BH}). This conjecture has been proved in several special cases, for example see (\cite{AD},\cite{HJY},\cite{D1},\cite{D2},\cite{DQ}) but it is still open in general.\\
\indent Let $I\subset S$ be a monomial ideal. It is well known that $\depth S/I \leq \depth S/\sqrt{I}$ (see the proof of \cite[Theorem 2.6]{HTT}) and equivalently $\depth I\leq \depth \sqrt{I}$. The first inequality holds also for $\sdepth$, that is $\sdepth S/I \leq \sdepth S/\sqrt{I}$ (see \cite [Theorem 1]{JP}). Moreover if $I\subset J$ are monomial ideals of $S$ then our Theorem \ref{0} says that $\sdepth J/I\leq \sdepth \sqrt{J}/\sqrt{I}$. In particular if $J/I$ is a kind of "Stanley-Cohen-Macaulay" module, that is $\sdepth J/I=\dim J/I$ then $\sqrt{J}/\sqrt{I}$ is too. The idea of the proof is inspired by \cite [Lemma 2.1]{DQ} .\\
\indent Next we give upper bounds for the intersection $I$ of two monomial primary ideals $Q$, $Q'$ (see Theorem \ref{16}). It is enough to find upper bounds for the intersection $J$ of two monomial prime ideals by Theorem \ref{0}. When these prime ideals have no common variables, this is done in Theorem \ref{7}, which says $\sdepth J\leq \frac{n+2}{2}$, using an idea of \cite[Lemma 2.2]{KSSY}  and the algorithm for computation of the Stanley depth given in \cite{HVZ}. For the general case we need to show that $\sdepth_{S[x_{n+1}]}(I,x_{n+1})\leq \sdepth_SI+1$, which is stated by the elementary Lemma \ref{10}. Using the lower bound given in \cite[Lemma 4.1]{DQ}  we noticed that our upper bound is reached when $Q$, $Q'$ are irreducible ideals, ht$(Q+Q')$ is odd and $\sqrt{Q}$, $\sqrt{Q'}$ have no common variables. In general our upper bound is big as shows Lemma \ref{11} and several examples.
\section{upper bounds for the stanley depth  }
Let $I$ and $J$ be two monomial ideals of $S$ such that $I\subset J$ and $\sqrt{I}$ and $\sqrt{J}$ be the radical ideals of $I$ and $J$ respectively. Let $G(I)$ be the minimal system of monomial generators of $I$. Then

\begin{Theorem}\label{0}$\sdepth(J/I)\leq \sdepth(\sqrt{J}/\sqrt{I}).$
\end{Theorem}
\begin{proof}
The ideals $I \ \hbox{and} \ J$ have an irredundant monomial decomposition  $J=\bigcap_{i=1}^{q}Q_i$ and $I=\bigcap_{j=1}^{r}Q'_j$, where $Q_i\hbox{'}s$ and $Q'_j\hbox{'}s$ are primary monomial ideals (if $J=S$ then $(Q_i)_i$ is the empty set of primary ideals). Let $T=K[y_1,\dots,y_n]$ and $a_i\in \NN$ be the maximum power of $x_i$ which appear in $\bigcup_{i=1}^{q}G(Q_i)\cup\bigcup_{j=1}^{r}G(Q_j)$, if no power of $x_i$ is there then take $a_i=0$. Set $a=\max_i(a_i)$. Let $P_i=(y_c:x_c\in \sqrt{Q_i}) \ \hbox{and} \ P'_j=(y_d:x_d\in \sqrt{Q'_j})$. These $P_i^{\hbox{,}}s$ and  ${P'}_{j}^{\hbox{,}}s$ are uniquely determined by $I \ \hbox{and} \ J$.  Let $\phi :T\longrightarrow S$ be the $K$-morphism given by $y_i\longrightarrow x_i^a$. Let $\mathcal{D} \ : J/I=\bigoplus _{i=1}^su_iK[Z_i]$ be a Stanley decomposition of $J/I$ such that $\sdepth(\mathcal{D})=\sdepth(J/I)$. Then we have \\$$L:=\phi ^{-1}(J)/\phi^{-1}(I)=\cap_{i=1}^{q}{P_i}/\cap_{j=1}^{r}{{P'}_j} \ ,$$ and $\phi$ defines an injection $\tilde{\phi}:L\longrightarrow J/I$. Thus  $$L=\tilde{\phi}^{-1}(\bigoplus_{i=1}^su_iK[Z_i])=\bigoplus_{i=1}^s \tilde{\phi}^{-1}(u_iK[Z_i]).$$ Suppose that $\tilde{\phi}^{-1}(u_jK[Zj])\neq 0$, then $u_j=x_{l_1}^{b_{l_1}}\dots \ x_{l_m}^{b_{l_m}}$ ,  $b_{l_k}\in \ZZ_+$ is such that if $x_{l_k}\notin Z_j$ then $a \ | \ b_{l_k}$ where $1\leq l_k \leq n$, let us say $b_{l_k}=ac_{l_k}$ for some $c_{l_k} \in \ZZ_+$. Denote $c_{l_k}=\lceil\frac{b_{l_k}}{a}\rceil$ when $x_{l_k}\in Z_j$. We get $$\tilde{\phi}^{-1}(u_jk[Z_j])=y_{l_1}^{c_{l_1}}\dots \ y_{l_m}^{c_{l_m}}k[V_j] \ ,$$ where $V_j=\{y_{l} : 1\leq l\leq n \ , \ x_{l}\in Z_j\}$. Thus $\tilde{\phi}^{-1}(u_jk[Z_j])$ is a Stanley space of $L$ and so $\mathcal{D}$ induces a Stanley decomposition $\mathcal{D'}$ of $L$ such that $$\sdepth(J/I)=\sdepth(\mathcal{D})\leq \sdepth(\mathcal{D'})\leq \sdepth(L).$$
Now let us define an isomorphism $\psi:T\longrightarrow S$ by $\psi(y_i)=x_i$. Under this isomorphism, we have $$L\cong \cap_{i=1}^{q}\sqrt{Q_i}/\cap_{j=1}^{r}\sqrt{Q'_j} \ .$$ Thus $$\sdepth(L)= \sdepth(\cap_{i=1}^{q}\sqrt{Q_i}/\cap_{j=1}^{r}\sqrt{Q'_j})=\sdepth(\sqrt{J}/\sqrt{I}).$$ Hence $$\sdepth(J/I)\leq \sdepth(\sqrt{J}/\sqrt{I}).$$
\end{proof}
\begin{Corollary}\label{1}
Let $I\subset S$ be a monomial ideal and $\sqrt{I}$ be its radical. Then $\sdepth(S/I)\leq \sdepth(S/\sqrt{I})$ and $\sdepth(I)\leq \sdepth(\sqrt{I})$.
\end{Corollary}
\noindent First part of the above corollary is already done in \cite[Theorem 1]{JP}.
\begin{Corollary}\label{2}
Let $I$ and $J$ be two monomial ideals of $S$ such that $I\subset J$ and $\sqrt{I}$ and $\sqrt{J}$ be the radical ideals of $I$ and $J$ respectively. If $\sdepth(J/I)=\dim(J/I)$. Then $\sdepth(\sqrt{J}/\sqrt{I})=\dim(\sqrt{J}/\sqrt{I})$.
\end{Corollary}
\begin{proof}
Given that  $\dim(J/I)=\sdepth(J/I)$,  and   from   Theorem \ref{0},  we have $$\dim(J/I)=\sdepth(J/I)\leq \sdepth(\sqrt{J}/\sqrt{I})\leq \dim(\sqrt{J}/\sqrt{I})$$ by \cite{SJ}. Also since $\dim(\sqrt{J}/\sqrt{I})=\dim(S/\sqrt{I}:\sqrt{J})\leq\dim(S/\sqrt{I:J})=\dim(S/I:J)=\dim(J/I)$, we get  $\sdepth(\sqrt{J}/\sqrt{I})=\dim(\sqrt{J}/\sqrt{I})$.
\end{proof}
\noindent The inequality given by Theorem \ref{0} can be strict as shows the following example.
\begin{Example}
{\em Let $I\subset K[x_1,x_2]$ and $I=(x_1^2,x_1x_2)$. Since $I=x_1^2K[x_1,x_2]\oplus x_1x_2K[x_2]$, we see that $\sdepth(I)=1$. Since $\sqrt{I}=(x_1)$, we have $\sdepth(\sqrt{I})=2$.
}
\end{Example}

A lower bound for Stanley depth of intersection of two monomial primary ideals is discussed in \cite{DQ}. We give an upper bound for the Stanley depth of intersection of two monomial primary ideals. Let $Q$ and $Q'$ be any two monomial primary ideal of $S$, then after renumbering the variables we can always assume that $\sqrt{Q}=(x_1,\dots,x_t)$ and $\sqrt{Q'}=(x_{r+1},\dots,x_p)$, where $0\leq r\leq t\leq p\leq n$. We start with the case $n=p$. A special case is given by the following:
\begin{Lemma}\label{4}
Let $Q$ and $Q'$ be two monomial primary ideals with $\sqrt{Q}=(x_1,\dots,x_t)$ and $\sqrt{Q'}=(x_{1},\dots,x_n)$. Then $$\sdepth(Q\cap Q')\leq n-\lfloor\frac{t}{2}\rfloor.$$
\end{Lemma}
\begin{proof}
$\sqrt{Q}\subseteq \sqrt{Q'}$ implies that $\sqrt{Q}\cap\sqrt{Q'}=\sqrt{Q}$ by \cite[Theorem 1.3]{C1} $\sdepth(\sqrt{Q})=n-\lfloor\frac{t}{2}\rfloor.$ So by Corollary \ref{1} it follows that $\sdepth(Q\cap Q')\leq n-\lfloor\frac{t}{2}\rfloor.$
\end{proof}

To find an upper bound for Stanley depth of $Q\cap Q'$ it is necessary to find an upper bound for $\sqrt{Q} \cap \sqrt{Q'}$. We consider two cases \\
\emph{\textbf{Case(1):}} $\sqrt{Q}=(x_1,\dots,x_t)$ and $\sqrt{Q'}=(x_{t+1},\dots,x_n)$ .\\
\emph{\textbf{Case(2):}} $\sqrt{Q}=(x_1,\dots,x_t)$ and $\sqrt{Q'}=(x_{r+1},\dots,x_n)$ where $1<r<t< n$ .\\

First we consider Case(1) where the proof idea comes from \cite{KSSY}. We recall the method of Herzog et al. \cite{HVZ} for computing the Stanley depth of a squarefree monomial ideal $I$ using posets. Let $G(I)=\{v_1,\dots,v_m\}$ be the set of minimal monomial generators of $I$. The characteristic poset of $I$ with respect to $g=(1,\dots,1)$  (see \cite{HVZ}), denoted by $\mathcal{P}_I^{(1,\dots,1)}$ is in fact the set $$\mathcal{P}_I^{(1,\dots,1)}=\{C\subset [n] \ | \ C  \ \text{contains supp($v_i$) for some $i$}\} \ ,$$ where  $\supp(v_i)=\{j:x_j|v_i\}\subseteq[n]:=\{1,\dots,n\}$. For every $A,B\in \mathcal{P}_I^{(1,\dots,1)}$ with $A\subseteq B$, define the interval $[A,B]$ to be $\{C\in \mathcal{P}_I^{(1,\dots,1)}:A\subseteq B\subseteq C\}$. Let $\mathcal{P}:\mathcal{P}_I^{(1,\dots,1)}=\cup_{i=1}^r[C_i,D_i]$ be a partition of $\mathcal{P}_I^{(1,\dots,1)}$, and for each $i$, let $c(i)\in \{0,1\}^n$ be the tuple such that $\supp(x^{c(i)})=C_i$. Then there is a Stanley decomposition $\mathcal{D}(\mathcal{P})$ of $I$ $$\mathcal{D}(\mathcal{P}):I=\bigoplus_{i=1}^sx^{c{(i)}}K[\{x_k|k\in D_i\}].$$ Clearly $\sdepth\mathcal{D}(\mathcal{P})=\min\{|D_1|,\dots,|D_s|\}$. It is shown in \cite{HVZ} that $$\sdepth(I)=\max\{\sdepth\mathcal{D}(P) \ | \ \mathcal{P} \text { is a partition of} \ \mathcal{P}_I^{(1,\dots,1)}\}.$$
\noindent An easy case which is enough when $n\leq 3$ is given by the following:
\begin{Lemma}\label{5}
Let $Q$ and $Q'$ be two monomial primary ideals with $\sqrt{Q}=(x_1)$ and $\sqrt{Q'}=(x_2,\dots,x_n)$. Then $$\sdepth(Q\cap Q')\leq1+\lceil\frac{n-1}{2}\rceil.$$
\end{Lemma}
\begin{proof}
Let $I$ be a monomial ideal and $v=GCD(u |u \in G(I))$. Then $I=vI'$ where $I'=(I:v)$. By \cite[Proposition  1.3(2)]{C2}  $\sdepth(I)=\sdepth(I')$. Since in our case  $v=x_1$ then we have $\sdepth(\sqrt{Q}\cap \sqrt{Q'})=\sdepth((\sqrt{Q}\cap \sqrt{Q'}) :x_1)=\sdepth(\sqrt{Q'})=n-(n-1)+\lceil\frac{n-1}{2}\rceil=1+\lceil\frac{n-1}{2}\rceil\ \  \hbox{by \ \cite[Theorem 1.3]{C1}}$. Hence from Corollary \ref{1} the result follows.
\end{proof}
\begin{Remark}\label{6}
{\em If in Lemma \ref{5} $Q$ and $Q'$ are irreducible monomial ideals then $\sdepth(Q\cap Q')=1+\lceil\frac{n-1}{2}\rceil$ by \cite[Theorem 1.3(2)]{C2}. Note that $\sdepth(Q\cap Q')=\frac{n}{2}+1$ if $n$ is even.}
\end{Remark}
\begin{Theorem}\label{7}
Let $Q$  and  $Q'$  be  two  primary  ideals  with $\sqrt{Q}=(x_1,\dots,x_t)$ and $\sqrt{Q'}=(x_{t+1},\dots,x_n)$, where $t\geq 2 \ and \ n\geq4$. Then
$$\sdepth(Q\cap Q')\leq \frac{n+2}{2} \ .$$
\end{Theorem}
\begin{proof}
\textbf{Case} $t=2 \ , \ n=4$. Applying Lemma \ref{5} it is enough to consider the case $\sqrt{Q}=(x_1,x_2)$ and $\sqrt{Q'}=(x_3,x_4)$. We have
\begin{multline*}
 \sqrt{Q}\cap \sqrt{Q'}=x_1x_3K[x_1,x_3,x_4]\oplus x_1x_4K[x_1,x_2,x_4]\oplus x_2x_3K[x_1,x_2,x_3]\\ \oplus x_2x_4K[x_2,x_3,x_4]\oplus  x_1x_2x_3x_4K[x_1,x_2,x_3,x_4].
\end{multline*}
The above Stanley decomposition shows that $\sdepth(\sqrt{Q}\cap \sqrt{Q'})=3$ because $\sqrt{Q}\cap \sqrt{Q'}$ is not principal. Then by Corollary \ref{1}  $\sdepth(Q\cap Q')\leq 3$, which is enough. For the remaining cases we proceed as follows:\\

Note that $\sqrt{Q}\cap \sqrt{Q'}$ is a squarefree monomial ideal generated in momomials of degree 2. Let $k:=\sdepth(\sqrt{Q}\cap \sqrt{Q'})$. The poset $P_{\sqrt{Q}\cap \sqrt{Q'}}$ has a partition $\mathcal{P}$ : $P_{\sqrt{Q}\cap \sqrt{Q'}}= \bigcup_{i=1}^{s}$ $[C_i,D_i]$ , satisfying $\sdepth(\mathcal{D}(\mathcal{P}))=k$ where $\mathcal{D(\mathcal{P})}$ is the Stanley decomposition of $\sqrt{Q}\cap \sqrt{Q'}$ with respect to partition $\mathcal{P}$. For each interval $[C_i,D_i]$ in $\mathcal{P}$ with $|C_i|=2$ we have $|D_i|\geq k$. There are $|D_i|-|C_i|$ subsets of cardinality $3$ in this interval. These intervals are disjoint.\\\textbf{Case} $t=2$ , $n\geq 5$\\ Since the number of subsets of cardinality 3 from that intervals $[C_i,D_i]$ with $|C_i|=2$ is at least
$$a:=\begin{bmatrix}
\begin{pmatrix}
  n \\
  2 \\
\end{pmatrix}
-\begin{pmatrix}
   n-2 \\
   2 \\
 \end{pmatrix}
-1
\end{bmatrix}
(k-2)$$ it follows that
$$a\leq
\begin{pmatrix}
  n \\
  3 \\
\end{pmatrix}
-\begin{pmatrix}
   n-2 \\
   3 \\
 \end{pmatrix}
.$$
\\
We get\\
\begin{equation*}
k\leq \frac{
\begin{pmatrix}
  n \\
  3 \\
\end{pmatrix}
-\begin{pmatrix}
   n-2 \\
   3\\
 \end{pmatrix}
}
{
\begin{pmatrix}
  n \\
  2 \\
\end{pmatrix}
-\begin{pmatrix}
   n-2 \\
   2 \\
 \end{pmatrix}
-1
}+2=\frac{n+2}{2} \ .
\end{equation*}
Thus we have $$\sdepth (Q\cap Q')\leq\sdepth(\sqrt{Q}\cap\sqrt{Q'})\leq\frac{n+2}{2} \ .$$\\
\noindent\textbf{Case} $t\geq 3 \ , \ n\geq 6$.\\ Now the number of subsets of cardinality 3 from the intervals $[C_i,D_i]$ with $|C_i|=2$ is at least
$$b:=\begin{bmatrix}
\begin{pmatrix}
  n \\
  2 \\
\end{pmatrix}
-\begin{pmatrix}
   t \\
   2 \\
 \end{pmatrix}
-\begin{pmatrix}
   n-t \\
   2 \\
 \end{pmatrix}  \\
\end{bmatrix}
(k-2)$$ and it follows $$
b\leq\
\begin{pmatrix}
  n \\
  3 \\
\end{pmatrix}
-\begin{pmatrix}
   t \\
   3 \\
 \end{pmatrix}
-\begin{pmatrix}
   n-t \\
   3 \\
 \end{pmatrix}.
$$
We get
\begin{equation*}
k\leq\frac{
\begin{pmatrix}
  n \\
  3 \\
\end{pmatrix}
-\begin{pmatrix}
   t \\
   3 \\
 \end{pmatrix}
-\begin{pmatrix}
   n-t \\
   3 \\
 \end{pmatrix}}
{
\begin{pmatrix}
  n \\
  2 \\
\end{pmatrix}
-\begin{pmatrix}
   t \\
   2 \\
 \end{pmatrix}
-\begin{pmatrix}
   n-t \\
   2 \\
 \end{pmatrix}
}+2=\frac{n+2}{2} \ .
\end{equation*}
Thus we have $$\sdepth(Q\cap Q')\leq\sdepth(\sqrt{Q}\cap \sqrt{Q'})\leq \frac{n+2}{2}.$$
\end{proof}
\begin{Corollary}\label{8}
Let $Q$ and $Q'$ be two irreducible monomial ideals such that $\sqrt{Q}=(x_1,\dots,x_t)$ and $\sqrt{Q'}=(x_{t+1},\dots,x_n)$. Suppose that $n$ is odd. Then $\sdepth(Q\cap Q')=\lceil\frac{n}{2}\rceil$.
\end{Corollary}
\begin{proof}
Using \cite[Lemma 4.1]{DQ}, we have, $\frac{n}{2}\leq \hbox{sdepth}(Q\cap Q')\leq\frac{n}{2}+1$ and so we get $\sdepth(Q\cap Q')=\lceil\frac{n}{2}\rceil$.
\end{proof}
\begin{Corollary}\label{9}
Let $Q$ and $Q'$ be two irreducible monomial ideals such that $\sqrt{Q}=(x_1,\dots,x_t)$ and $\sqrt{Q'}=(x_{t+1},\dots,x_n)$. Suppose that $n$ is even. Then
$$\sdepth(Q\cap Q')=
\left\{
  \begin{array}{ll}
    \frac{n}{2}+1, &  if \ \ \ t \ is \ odd; \\
     & \hbox{} \\
    \frac{n}{2}\ or\ \frac{n}{2}+1, & \hbox{if} \ \ \ t \ is \ even.
  \end{array}
\right.$$
\end{Corollary}
\begin{proof}
Using again \cite[Lemma 4.1]{DQ}, we have, $\lceil\frac{t}{2}\rceil+\lceil\frac{n-t}{2}\rceil\leq \sdepth(Q\cap Q')\leq\frac{n}{2}+1$.\\\textbf{\emph{Case}} $t$ is odd ,\\If  $t$ is odd then $\lceil\frac{t}{2}\rceil+\lceil\frac{n-t}{2}\rceil=\frac{n}{2}+1$ and so we get $\sdepth(Q\cap Q')=\frac{n}{2}+1$ in this case.\\
\emph{\textbf{Case}} $t$ is even,\\If $t$ is even then $\lceil\frac{t}{2}\rceil+\lceil\frac{n-t}{2}\rceil=\frac{n}{2}$  and we get $\sdepth(Q\cap Q')=\frac{n}{2} \ or \ \frac{n}{2}+1.$
\end{proof}
\noindent We need next the following elementary lemma:
\begin{Lemma}\label{10}
Let $I$ be a monomial ideal of $S$, and let $I'=(I,x_{n+1})$ be the monomial ideal of $S'=S[x_{n+1}]$. Then
$$\sdepth_S(I)\leq\sdepth_{S'}(I')\leq \sdepth_S(I)+1.$$
\end{Lemma}
\begin{proof}
Since $I'=(I,x_{n+1})$, so $I'\cap S=I$. Now let $\mathcal{D} \ : \ I'=\bigoplus_{i=1}^{r} u_iK[Z_i]$ be the Stanley decomposition of $I'$, with $\sdepth(\mathcal{D})=\sdepth_{S'}(I').$
Since
$$I=I'\cap S=\bigoplus_{i=1}^{r} u_iK[Z_i]\cap S$$
$$ \ \ \ \  \ \ \ \ \ \ \ \ \ \ \ \ \ \ \ \ \ \ \ \ \ \ \ \ \ \ \ \ =\bigoplus_{x_{n+1} \notin supp(u_i) } u_i K[Z_i \ \backslash \ \{x_{n+1}\}]$$
we conclude that  $$\sdepth_{S'}(I')\leq \sdepth_S(I)+1.$$
For the other inequality note that a Stanley decomposition $\mathcal{D}:I=\bigoplus_iu_iK[Z_i]$ with $\sdepth\mathcal{D}=\sdepth_S(I)$ induces a Stanley decomposition\\$$\mathcal{D'}:I'=\bigoplus_iu_iK[Z_i]\bigoplus x_{n+1}S'$$ with $\sdepth\mathcal{D'}=\sdepth(I).$
\end{proof}
\begin{Example}
{\em
If $I$ is any monomial complete intersection ideal of $S$ with $|G(I)|=m$ then   $\sdepth_S(I)=n-\lfloor\frac{m}{2}\rfloor$ by \cite{Y}. Note that $I'=(I,x_{n+1})$ is again a monomial complete intersection ideal of $S'$, so we have $\sdepth_{S'}(I')=n+1-\lfloor\frac{m+1}{2}\rfloor$. Now if $m$ is odd then $\sdepth_{S'}(I')=\sdepth_S(I)$ and  $\sdepth_{S'}(I')=\sdepth_S(I)+1$ if $m$ is even.
}
\end{Example}
\begin{Proposition}\label{12}
Let $Q$ and $Q'$ be two primary ideals with $\sqrt{Q}=(x_1,\dots,x_t)$ and $\sqrt{Q'}=(x_{r+1},\dots,x_n)$, where $1<r\leq t<n$,  $n\geq 4$. Then $$\sdepth(Q\cap Q')\leq \frac{n+t-r+2}{2} \ .$$
\end{Proposition}
\begin{proof}
Let $S'=K[x_1,\dots,x_{r},x_{t+1},\dots,x_n]$ and $q=Q\cap S'$ , $q'=Q'\cap S'$. Then $$\sdepth(q\cap q')\leq \frac{(n-t+r)+2}{2} \ \ \ \ \text{by Theorem \ref{7}}.$$
But $\sdepth(Q\cap Q')\leq \sdepth(q\cap q')+t-r$ applying Lemma \ref{10} by recurrence and the inequality follows.
\end{proof}
The bound given by Proposition \ref{12} sounds reasonable for $t=r$, otherwise seems to be too big as shows our Example \ref{13} and Lemma \ref{11}.
\begin{Example}\label{13}
{\em Let $\mathfrak{m}=(x_1,\dots,x_n)$ be the maximal ideal of $S$. Then $\sdepth(\mathfrak{m})=\lceil\frac{n}{2}\rceil$. Now let $S'=S[x_{n+1},\dots,x_{n+r}]$ and $\mathfrak{m'}=(\mathfrak{m},x_{n+1},\dots,x_{n+r}).$ We see that $\sdepth_{S'}(\mathfrak{m'})=\lceil\frac{n+r}{2}\rceil$, which is much smaller than $\sdepth_S\mathfrak{m}+r$.
}
\end{Example}
\begin{Lemma}\label{11}
Let $Q$ and $Q'$ be two primary monomial ideals with $\sqrt{Q}=(x_1,\dots,x_{n-1})$ and  $\sqrt{Q'}=(x_2,\dots,x_n)$. Then  $$\sdepth(Q\cap Q')\leq n-\lfloor\frac{n-1}{2}\rfloor.$$
\end{Lemma}
\begin{proof}
Since  $\sqrt{Q}\cap \sqrt{Q'}=(x_1x_n,x_2,\dots,x_{n-1})$ is a complete intersection ideal we have $\sdepth(\sqrt{Q}\cap \sqrt{Q'})=n-\lfloor\frac{n-1}{2}\rfloor$ by \cite[Theorem 2.4]{Y}. So by Corollary \ref{1} we have  that $$\sdepth(Q\cap Q')\leq n-\lfloor\frac{n-1}{2}\rfloor.$$
\end{proof}
\noindent Another possible bound is given below.
\begin{Proposition}\label{14}
Let $Q$ and $Q'$ be two monomial primary ideals with $\sqrt{Q}=(x_1,\dots,x_t)$ and $\sqrt{Q'}=(x_{r+1},\dots,x_n)$ where $1<r\leq t<n$. Then $$\sdepth(Q\cap Q')\leq \min\{n-\lfloor\frac{t}{2}\rfloor, \ n-\lfloor\frac{n-t}{2}\rfloor\}.$$
\end{Proposition}
\begin{proof}
Let $I'=(x_1,\dots,x_t)\cap(x_r,\dots,x_n)$ be an ideal of $S$. Let $\phi : S\longrightarrow S'$ where $S'=K[x_1,\dots,x_{n-1}]$ be the homomorphism given by  $\phi(x_i) = x_i$ for $i \leq n-1$ and $\phi(x_n) = 1$. We see that $I=\phi(I')=(x_1,\dots,x_t)$ where $I\subset S'$. Then by \cite[Lemma 2.2]{C1}, we have
$$\sdepth_S(I')\leq \sdepth_{S'}(I)+1.$$
Since $$\sdepth_{S'}(I)=(n-1)-t+\lceil\frac{t}{2}\rceil$$ by \cite[Theorem 1.3]{C1} we have $$\sdepth_S(I')\leq (n-1)-t+\lceil\frac{t}{2}\rceil+1=n-\lfloor\frac{t}{2}\rfloor.$$ By Corollary \ref{1} we have$$\sdepth(Q\cap Q')\leq \sdepth_S(I')\leq n-\lfloor\frac{t}{2}\rfloor$$ and similarly $$\sdepth(Q\cap Q')\leq n-\lfloor\frac{n-t}{2}\rfloor.$$
\end{proof}
\begin{Example}\label{17}
{\em Let $I=(x_1,\dots,x_6)\cap(x_3,\dots,x_8)$, $n=8$. Then by Proposition \ref{14}, we have   $$\sdepth(I)\leq 8-\lfloor\frac{6}{2}\rfloor=5.$$ But using Proposition \ref{12}, we have $\sdepth(I)\leq 7$, which shows that Proposition \ref{14} gives better  bound than Proposition \ref{12} in this case.}
\end{Example}
\begin{Example}
{\em Let $I=(x_1,\dots,x_5)\cap (x_5,\dots,x_9)$,  $n=9.$  Then by Proposition \ref{14}, we have $\sdepth(I)\leq7$. But using Proposition \ref{12}, we have $\sdepth(I)\leq 6.$
}
\end{Example}

These examples show somehow that the upper bound given by Proposition \ref{12} is good if less number of variables are in common. The upper bound of Proposition \ref{14} seems to be better if we have large number of variables in common.
\begin{Theorem}\label{16}
Let $Q$ and $Q'$ be two monomial primary ideals with $\sqrt{Q}=(x_1,\dots,x_t)$ and $\sqrt{Q'}=(x_{r+1},\dots,x_p)$, where $1<r\leq t<p\leq n$, $n\geq 4$. Then $$\sdepth(Q\cap Q')\leq \min\{\frac{2n+t-p-r+2}{2},\  n- \lfloor\frac{t}{2}\rfloor, \ n-\lfloor\frac{p-t}{2}\rfloor\}.$$ The inequality becomes equality if $t=r$, $n$ is odd and $Q$, $Q'$ are irreducible.
\end{Theorem}
\begin{proof}
Let $S'=K[x_1,\dots,x_p]$ and $q=Q\cap S'$, $q'=Q'\cap S'$. Then $$\sdepth(q\cap q')\leq \frac{p-t+r+2}{2} \ \ \  \text{ by Proposition \ref{12}}$$  and  $$\sdepth(q\cap q')\leq \min\{p-\lfloor\frac{t}{2}\rfloor, \ p-\lfloor\frac{p-t}{2}\rfloor\} \ \ \text{ by Proposition \ref{14}}.$$ Thus $$\sdepth(Q\cap Q')=\sdepth(q\cap q')+n-p\leq \frac{2n+t-p-r+2}{2}$$ and $$\sdepth(Q\cap Q')=\sdepth(q\cap q')+n-p\leq \min\{n-\lfloor\frac{t}{2}\rfloor, \ n-\lfloor\frac{p-t}{2}\rfloor\}$$ by \cite[Lemma 3.6]{HVZ}. For the second statement apply Corollary \ref{8} and \cite[Lemma 3.6]{HVZ}.
\end{proof}

\begin{Example}
{\em Let $I=(x_1^{a_1},x_2^{a_2})\cap (x_2^{a_2},x_3^{a_3},\dots,x_n^{a_n})$, where $n=2k$ and $k\geq2$. Then $\sdepth(I)=k+1$. Indeed, from Proposition \ref{12} we have $$\sdepth(I)\leq\lfloor\frac{2k+3}{2}\rfloor= k+1.$$ Also we know from \cite{O} that $$\sdepth(I)\geq n-\lfloor\frac{|G(I)|}{2}\rfloor=2k-\lfloor\frac{2k-1}{2}\rfloor=2k-k+1=k+1,$$ which is enough.}
\end{Example}

\end{document}